\documentclass[11pt]{amsart}

\usepackage{graphicx}
\usepackage{amscd}
\usepackage{amsmath}
\usepackage{amsfonts}
\usepackage{amssymb}
\usepackage{amsxtra}
\usepackage{xspace}
\usepackage{array}

\theoremstyle{plain}
\newtheorem{theorem}{Theorem}
\newtheorem{lemma}{Lemma}
\newtheorem{corollary}{Corollary}
\newtheorem{proposition}{Proposition}

\newtheorem{remark}{Remark}

\theoremstyle{definition}
\newtheorem{definition}{Definition}

\numberwithin{equation}{section}

\newcommand{\be}{\begin{enumerate}}
\newcommand{\ee}{\end{enumerate}}
\newcommand{\beq}{\begin{equation}}
\newcommand{\eeq}{\end{equation}}
\newcommand{\bprop}{\begin{proposition}}
\newcommand{\eprop}{\end{proposition}}

\newcommand{\g}[2][W^J]{\gamma(#1, #2)} 
\newcommand{\rs}[1][]{\Phi^{#1}}
\newcommand{\lc}[1][W\p \,]{[#1]}  

\newcommand{\p}{^\prime}
\newcommand{\wth}{W_{(3)}}
\newcommand{\sth}{S_{(3)}}
\newcommand{\tth}{T_{(3)}}
\newcommand{\vth}{V_{(3)}}
\newcommand{\wj}{W^J}
\newcommand{\ovw}{{\widetilde{W}\p}}
\newcommand{\oh}{\mathcal{O}_1}

\newcommand{\complex}{\mathbb{C}}
\newcommand{\reals}{\mathbb{R}}
\newcommand{\nat}{\mathbb{Z}^{\geq 0}}
\newcommand{\integers}{\mathbb{Z}}

\newcommand{\pfbegin}{\noindent {\em Proof:} }

\begin{document}

\title[]{On growth types of quotients of Coxeter groups by parabolic subgroups}
\author{Sankaran Viswanath}
\address{Department of Mathematics\\
University of California\\
Davis, CA 95616, USA}
\email{svis@math.ucdavis.edu}

\begin{abstract}
The principal objects studied in this note are Coxeter groups $W$ that are
 neither finite nor affine. A well known result of de la Harpe asserts
 that such groups have exponential growth. We consider quotients of $W$ by its 
parabolic subgroups and by a certain class of reflection subgroups. We
show that these quotients have exponential growth as well. To achieve
 this, we use a theorem of Dyer to construct a reflection subgroup of
 $W$ that is  isomorphic to the universal Coxeter group on three generators. 
The results are all proved under the restriction that
the Coxeter diagram of $W$  is simply laced, and some remarks
made on how this restriction may be relaxed.
\end{abstract}
\maketitle

\section{Introduction}
Let $W$ be a finitely generated group and $S$ be a finite set of
generators for $W$. The {\em growth function} $\gamma(m)$ is the number of elements of
$W$ expressible as a word of length $m$ or less in $S \cup S^{-1}$. We
say that $W$ has (i) polynomial growth if $\exists \, C \in \reals^{>0}$
and $d \in \nat$ such that $\gamma(m) \leq Cm^d\, \forall m \geq 0$,
(ii) exponential growth if $\exists \, \lambda >1$ such that
$\gamma(m) \geq \lambda^m\, \forall m \geq 0$ and (iii) intermediate
growth otherwise. For  preliminaries on growth types of finitely
generated groups, we refer the reader to de la Harpe's monograph
\cite[Chap. VI, VII]{pdlhbk} or to section \ref{growth} below.


We now specialize to the case where $(W,S)$ is an irreducible Coxeter
system. If $W$ is a finite or affine Coxeter group, it can be easily seen to
have polynomial growth. When $W$ is an infinite, non-affine Coxeter
group, it is a classical result of de la Harpe \cite{pdlhpap} that $W$
has exponential growth.
In this note, we consider the latter case. We will prove the slightly stronger result that for any
proper parabolic subgroup $W_J$ of $W$, the quotient $W/W_J$  has
exponential growth too.
This quotient can be identified with the set $W^J$ of minimal length
left coset representatives; by the ``growth type of $W/W_J\,$'' we will
mean the growth type of the subset $W^J$ of $W$.

We remark that this assertion about the quotient $W/W_J$ does not follow directly from
the exponential growth of $W$ given by de la Harpe's theorem. The
group $W$
could have parabolic subgroups $W_J$ that are infinite, non-affine and thus of
 exponential growth themselves. For such $W_J$, the growth type of the quotient
 $W/W_J$ is not apriori determined.

 The main ingredient in our approach to the growth of $W/W_J$ is the
 work of Deodhar \cite{deodhar2} and Dyer \cite{dyer1} on reflection subgroups
 of Coxeter groups. Using a criterion of Dyer, we construct a specific
 reflection subgroup of $W$; this subgroup will turn out to have two
 properties of interest to us : (i) it is isomorphic to the ``universal''
Coxeter group on three generators, and (ii) distinct reflections in this
 subgroup belong to distinct cosets in $W/W_J$.
These properties will enable us to deduce the exponential growth of $W/W_J$.

We apply this result on the growth of $W/W_J$ to study the growth of
more general quotients $W/W\p$, where $W\p$ is a reflection subgroup
of $W$. We identify a class of reflection subgroups $W\p$ of $W$
 for which the quotients $W/W\p$ have exponential growth. 

 We'll work throughout under the hypothesis that $W$ is simply laced; this
 restriction can however be relaxed and we indicate this in the
 relevant places (see remark \ref{light}).

Here's a quick outline of the rest of the article: in section 2, after
preliminaries on growth types, we state our main theorem. Section 3 is
concerned with the construction of the special reflection subgroup
mentioned above, and section 4 collects together some well known facts
about universal Coxeter groups. These facts are then applied to our
reflection subgroup to complete the proof of the main theorem in
section 5. In the final section, we study the more general quotients
$W/W\p$, where $W\p$ is a reflection subgroup satisfying some
additional hypothesis.

\section{Growth types} \label{growth}
\subsection{} 
We follow \cite[Chapter VI.C]{pdlhbk} :
\begin{definition}
Given a non decreasing sequence $(a_k)_{k\geq 0}$ of natural numbers,
its {\em exponential growth rate} is defined to be $\omega :=
\displaystyle\limsup_{k \to \infty} a_k^{1/k}$.
\end{definition}
Now suppose  $(W,S)$ is a Coxeter system and $F \subset W$ is a
subset such that $1 \in F$. 
For $k \geq 0$, set $\g[F]{k}:=\# \{f \in F: \ell(f) \leq k\}$ and 
 $\omega(F) := \displaystyle\limsup_{k \to \infty} \g[F]{k}^{1/k}$. Since $1
\leq \g[F]{k} \leq (\# S + 1)^k \; \forall k$, we have $1 \leq
\omega(F) \leq (\# S +1)$.

\begin{definition}\label{defexpgr}
We say that $F$ has exponential growth if $\omega(F) > 1$ and
subexponential growth otherwise.
\end{definition}
A special case of subexponential growth is {\em polynomial growth},
which occurs if $\exists\, C \in \reals^{>0}$ and $d \in \nat$ such that $\g[F]{k}
\leq C k^d$ for all $k \geq 0$. If $F$ is of subexponential growth and
not of polynomial growth, we say it has {\em intermediate growth}.

When $F = W$, the function $\g[W]{k}$ is submultiplicative i.e,
satisfies $\g[W]{k+l} \leq \g[W]{k}\cdot \g[W]{l}$. This implies (see
\cite[VI.56]{pdlhbk}) that $\lim_{k \to \infty} \g[W]{k}^{1/k}$ exists
and equals $\inf_{k \geq 0} \g[W]{k}^{1/k}$. Thus we get an equivalent
formulation: $W$ has exponential growth iff  $\exists\, \lambda > 1$ such that
$\g[W]{k} \geq \lambda^k$ for all $k \geq 0$. If $F$ is a proper
subset of $W$, then submultiplicativity need not hold, and we will be
content with definition \ref{defexpgr} for our notion of exponential growth.

\subsection{Rational generating functions and growth}
Given $\{1\} \subset F \subset W$ as above, let $\gamma_F(q) \in \complex[[q]]$ be the
generating function:
$$\gamma_F(q) := \sum_{k=0}^{\infty} \g[F]{k} q^k$$
Observe that $\omega(F)^{-1}$ is the radius of convergence of this
power series. For a Coxeter group $W$, there are many
natural choices of $F$ (e.g parabolic subgroups, their minimal coset
representatives) for which $\gamma_F(q)$ is a rational function. When this
happens, one clearly also has:
\begin{proposition}\label{rfgrowth}  Suppose $\gamma_F(q)$ is a rational function. Then 
$F$ has exponential growth iff $\gamma_F(q)$ has a pole $\xi$ with $0 <
  |\xi| < 1$.
\end{proposition}
See \cite[proposition 3.3]{wagreich} for the situation when $F$ has
polynomial growth.

\subsection{}
Let $(W,S)$ be an irreducible Coxeter system.
Let $W_J, J \subsetneq S$ be a parabolic subgroup of $W$, with
$\wj$ being the set of minimal length elements in left cosets of
$W_J$. Recall that each $w \in W$ can be uniquely written as $w =
\sigma \tau, \sigma \in \wj, \tau \in W_J$ with $\ell(w) = \ell(\sigma)
+ \ell(\tau)$. Our objective is to study the growth type of the
subset $\wj$.

When $W$ is finite or affine, it is easy to see that the set $W^J$ has the same growth type
as $W$. We consider the case where $W$ is an infinite, non-affine Coxeter
group. We will further assume that the
Coxeter diagram of $W$ is connected and simply laced i.e for each pair
$s \neq  s^{\prime} \in S$, $ss^{\prime}$ has order 2 or 3 in $W$.  
Our main theorem is the following:

\begin{theorem}\label{mainthm} Let $(W,S)$ be an irreducible Coxeter
  system. Suppose $W$ is infinite, non-affine and has a 
simply laced Coxeter diagram. Then for all $J \subsetneq S$, 
$\wj$ has exponential growth.
\end{theorem}
\begin{corollary}\label{cor-one}
Under the assumptions of theorem \ref{mainthm}, the Poincar\'{e} series (length
generating function) of $W^J$ has a pole $\xi$ with $0 < |\xi| < 1$.
\end{corollary}

\noindent
{\em Proof of Corollary \ref{cor-one}:}
If $W^J(q) = \sum_{\sigma \in W^J} q^{\ell(\sigma)}$ is the Poincar\'{e}
series of $W^J$, we have $\gamma_{W^J}(q) = W^J(q)/(1-q)$. The
corollary now follows from theorem \ref{mainthm} and proposition
\ref{rfgrowth}. \qed

\vspace{0.2cm}
Observe that if $W_J$ is a finite group, the assertion of Theorem
\ref{mainthm} is trivial. So we may as well assume that $W_J$ is
infinite. To show that $W/W_J$ has exponential growth, we will 
do two things: (A) construct a large (exponential in $m$) number of elements
in $W$ of length $\leq m$ and (B) show that these elements lie in
distinct left cosets of $W_J$.

To achieve step (B), we will employ the following nice result due to
Deodhar \cite{deodhar1}:
\begin{theorem}(Deodhar)\label{deo-thm1}
Let $(W,S)$ be a Coxeter system, $T:=\bigcup_{w \in W}
wSw^{-1}$ be the set of reflections and $J \subset S$. If $t_1,
t_2 \in T\backslash W_J$ with $t_1\neq t_2$,  then $t_1 W_J \neq t_2W_J$; i.e
distinct elements of $T\backslash W_J$ are  in  distinct left cosets of $W_J$.
\end{theorem}

The next proposition makes step (A) above more precise:
\begin{proposition}\label{mainprop}
Assume notation as in the statement of theorem~\ref{mainthm}. Suppose
also that $W_J$ is an infinite group. Then there exists a natural
number $M$ such that for all $k \geq 0$, $\exists\,$ at least $2^k$
elements $t \in T\backslash W_J$ with $\ell(t) \leq M(2k+1)$.
\end{proposition}
Given the truth of this proposition, we now have:

\noindent
{\em Proof of theorem~\ref{mainthm}:} For $t \in T$, let $[t] \in \wj$ denote the
unique minimal length element in $tW_J$. For $t \in T\backslash W_J$   as in proposition
\ref{mainprop}, $\ell([t]) \leq \ell(t)
\leq M(2k+1)$. Invoking theorem~\ref{deo-thm1}, we conclude  that there exist at
least $2^k$ elements $w \in \wj$ with $\ell(w) \leq M(2k+1)$ i.e
$\g{M(2k+1)} \geq 2^k$. This gives for $k \geq 1$:
$$\g{M(2k+1)}^{1/M(2k+1)} \geq 2^{k/M(2k+1)} \geq 2^{1/3M}$$
Thus $\omega(\wj) = \limsup_{k \to \infty} \g{k}^{1/k} \geq 2^{1/3M} >
1$. So $\wj$ has exponential growth. This completes the proof of
theorem \ref{mainthm}\qed

The next three sections will be devoted to a proof of
proposition~\ref{mainprop}.

\section{A reflection subgroup isomorphic to $\wth$} \label{reflnsubgp}
\subsection{}
As a first step toward proving proposition~\ref{mainprop}, we will
construct a reflection subgroup of $W$ that is isomorphic to the
universal Coxeter group $\wth = \langle  s_1, s_2, s_3: s_i^2=1
\,\forall i \rangle$ \cite{dyer2}
with Coxeter diagram  
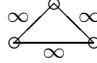
\setlength{\unitlength}{15pt}
$$\begin{picture}(2,2)(0,-.5)
\put(0,0){\circle{.25}}
\put(2,0){\circle{.25}}
\put(1,1){\circle{.25}}

\put(0,0){\line(1,0){2}}
\put(0,0){\line(1,1){1}}
\put(2,0){\line(-1,1){1}}

\put(1,-.25){\makebox(0,0){\scriptsize $\infty$}}
\put(.1,.65){\makebox(0,0){\scriptsize $\infty$}}
\put(1.9,.65){\makebox(0,0){\scriptsize $\infty$}}
\end{picture}$$
We collect together  the relevant facts about reflection subgroups from
Deodhar \cite{deodhar2}\ and Dyer \cite{dyer1}. We recall that the
elements of the set $T:=\bigcup_{w \in W} wSw^{-1}$ are called {\em
  reflections}.
\begin{definition}
Let $(W,S)$ be a Coxeter system.
A subgroup $W\p$ of $W$ generated by reflections is called a reflection subgroup.
\end{definition}
Reflection subgroups of Coxeter groups turn out to be Coxeter groups
in their own right. Specifically:
\begin{theorem}\label{deo-dyer} 
(Deodhar \cite{deodhar2}, Dyer \cite{dyer1}) Let $W\p$ be a reflection
  subgroup of $W$. Then $\exists\, S\p \subset W\p \cap T$ such that
  $S\p$ forms a set of Coxeter generators for $W\p$.
\end{theorem}

\subsection{}
Let $(W,S)$ be a Coxeter system with simply laced Coxeter diagram $X$.
We assume that the nodes of $X$ are labelled by
the elements of $S$. We let $V$
denote the geometric representation of $W$ \cite[\S 5.3]{humphreys}; 
$V$ has a basis $\{\alpha_s: s \in S\}$ and a
symmetric bilinear form $(,)$ determined by the conditions:
(i) $(\alpha_s, \alpha_s)=1 \; \forall s \in S$; (ii) $ (\alpha_{p}, \alpha_{q})=-1/2$
when $p \neq q \in S$ and the nodes $p$ and $q$ are connected by an edge
in $X$ and $(\alpha_p, \alpha_q)=0$ otherwise. The $W$ action on $V$ 
preserves the form $(,)$ and is determined by $s(v) = v - 2(v,\alpha_s)
\alpha_s$ where $s \in S, v \in V$.

Let $\rs(W) = W.\{\alpha_s: s \in S\} \subset V$ be the root
system of $W$ \cite[\S 5.4]{humphreys} and let $\rs[+](W)$
(resp. $\rs[-](W)$) denote the set of positive (resp. negative) roots
of W. The set $\rs[+](W)$ is in 1-1 correspondence with the set $T$
via $\alpha \mapsto s_{\alpha}$  defined by 
$s_{\alpha}(v):=v-2(v,\alpha) \alpha$. The following theorem of
Dyer \cite[Theorem 4.4]{dyer1} will be important in what follows:
\begin{theorem}\label{dyer}(Dyer)
Let $F \subset \rs[+](W)$, $S\p :=\{s_\alpha:
\alpha \in F\}$ and $W\p$ be the reflection subgroup of $W$ generated
by $S\p$. Suppose for all $\alpha \neq \beta \in F$, $(\alpha, \beta)
\in \{-\cos(\pi/n): n \in \mathbb{N}, n \geq 2\} \cup (-\infty, -1]$, then $S\p$ is precisely
    the set of Coxeter generators of $W\p$ given by theorem \ref{deo-dyer}.
\end{theorem}
\begin{corollary}\label{w3constr}
 If $\beta_i \,(i=1,2,3) \in \rs[+](W)$ are such that $(\beta_i, \beta_j)
 \leq -1$ for all $i \neq j$, then the reflection subgroup 
$W\p = \langle  s_{\beta_1}, s_{\beta_2}, s_{\beta_3}\rangle \subset W$ is
 isomorphic to $\wth$.
\end{corollary}
\noindent
\pfbegin By Theorem~\ref{dyer} the $s_{\beta_i}$ are the Coxeter
generators of $W\p$. It is an easy fact (see for e.g
 \cite[\S 5.3]{humphreys}) that if $(\beta_i, \beta_j) \leq -1$, then
$order(s_{\beta_i}s_{\beta_j}) = \infty$. Thus $W\p \cong \wth$. \qed

\subsection{}\label{thebetas}
We now assume the notation as in the statement of
Proposition~\ref{mainprop}. So (i) $(W,S)$ is an irreducible Coxeter
system, (ii) $W$ is infinite, non-affine with simply laced
Coxeter diagram $X$ and (iii) $J \subsetneq S$ with
$W_J$ infinite. We will use corollary~\ref{w3constr} to construct a
reflection subgroup of $W$ isomorphic to $\wth$.

In what follows, we will identify (without explicit mention) subsets
$K$ of $S$ with the corresponding subdiagram of $X$ formed by taking only the nodes
labelled by $K$, together with all edges between these nodes.

First, we decompose $J = \bigsqcup J_j$ where the $J_j$ are the {\em connected
 components} of $J$. Since $\prod_j W_{J_j} \cong W_J$ is
infinite, $\exists\, i$ such that $W_{J_i}$ is infinite.
Let $Z:=J_i$\,; $Z$  is thus not a diagram of {\em finite
 type}. It is a classical result (verifiable by hand)
that any connected, simply laced diagram either contains or is
contained in one of the affine simply laced diagrams $\widetilde{A}_n, n
\geq 2$, $\widetilde{D}_n, n \geq 4$, $\widetilde{E}_n, n=6,7,8$ (this
result can in fact be used to quickly classify the finite simply laced
Coxeter groups). Applying this 
to $Z$, one concludes that $Z$ must {\em contain} an affine diagram
 $Y$; this is because if $Z$ were properly contained in an affine
 diagram, then $Z$ would end up being of finite type.

Now pick $p \in S\backslash J$; clearly $p \not\in Y$.
 Since the Coxeter diagram 
$X$ is connected, we can pick a shortest path in $X$ between $p$ and
$Y$; i.e, $\exists\, s_o, s_1, \cdots, s_r \in S$ such that 
\be
\item $s_0=p$.
\item $s_r \in Y$.
\item $s_i$ and $s_{i+1}$ are connected by an edge in $X \; \forall i$.
\item $r$ is the smallest such natural number.
\ee
The minimality of $r$ is easily seen to imply the following:
\be
\item The $s_i, 0 \leq i \leq r$ are pairwise distinct.
\item $s_i \not\in Y$ for $i<r$.
\item $s_i$ and $s_j$ are not connected by an edge in $X$ if
  $|i-j|>1$.
\ee
Thus the subdiagram of $X$ formed by the nodes labelled $s_i$ is just the
  classical diagram $A_{r+1}$.

Now, if $K$ is a subset of $S$, we will naturally identify
$\rs(W_K)$ with the subset $W_K.\{\alpha_k: k \in K\}$ of $\rs(W)$. 
We recall that since $Y$ is an affine
diagram, there exists $\delta_Y \in \rs[+](W_Y)$ 
such that $(\delta_Y, \alpha_q)=0 \; \forall \text{ nodes } q \in
Y$. Further, if $\delta_Y = \sum_{q \in Y} c_q \alpha_q$, then we have
$c_q \geq 1$ for all $q \in Y$.

Now, define positive roots $\beta_i \, (i=1,2,3) \in \rs[+](W)$ as
follows:
\begin{align*}
\beta_1 &:=\sum_{i=0}^{r-1} \alpha_{s_i} = s_{r-1}\cdots
s_{2}s_{1}(\alpha_{p})\\
\beta_2 &:= \alpha_{s_r} + \delta_Y \\
\beta_3 &:= -\alpha_{s_r} + 3\delta_Y 
\end{align*}
Observe that the $\beta_i$ are linearly independent.
By the well known characterization of positive roots of affine Coxeter
groups, we have  $\beta_2, \beta_3 \in \rs[+](W_Y) \subset
\rs[+](W)$. Further $\beta_1$ has been explicitly demonstrated to be an
element of $\rs[+](W)$. For $i=2,3$, if we write $\beta_i = \sum_{q \in Y} 
 c_q^{(i)} \alpha_q$, then we have $c_{s_r}^{(i)} \geq 2$ for both
 values of $i$. So 
\begin{align*}
(\beta_1, \beta_2) &= (\sum_{i=0}^{r-1} \alpha_{s_i}, \alpha_{s_r} +
 \delta_Y ) = (\sum_{i=0}^{r-1} \alpha_{s_i}, c_{s_r}^{(2)}\alpha_{s_r} + \sum_{q \in Y, q\neq
 s_r}  c_q^{(2)} \alpha_q)\\
&\leq (\alpha_{s_{r-1}}, 2\alpha_{s_r}) = -1
\end{align*}
Similarly $(\beta_1, \beta_3) \leq -1$ too. Finally $(\beta_2, \beta_3)
= (\alpha_{s_r} + \delta_Y, -\alpha_{s_r} +3 \delta_Y) =
(\alpha_{s_r}, -\alpha_{s_r}) = -1$.
Corollary~\ref{w3constr} can now be applied to deduce:
\bprop\label{sbetas}
$\langle s_{\beta_1}, s_{\beta_2},  s_{\beta_3}\rangle \cong \wth $
\eprop
\begin{remark}\label{light}
Suppose the Coxeter diagram of $W$ is not simply laced, but $J$ contains
an affine subdiagram $Y$, then it is clear that the above construction
can still be carried out with some simple modifications. In particular
if for all pairs $s \neq s\p \in S, (ss\p)^{m_{ss\p}} =1$ with $m_{ss\p}=2, 3$
  or $\infty$ and rank $W \geq 3$, the above construction works for
  all $ J \subsetneq  S$ with $\# W_J =\infty$.
\end{remark}
 We also note the following interesting corollary to the above construction:
\begin{corollary} \label{expo}
Let $(W,S)$ be an irreducible Coxeter system. 
If $W$ is  infinite, non-affine and  simply laced, then $W$
contains a reflection subgroup isomorphic to the universal Coxeter
group $\wth$.
\end{corollary}

\noindent
\pfbegin Let $X$ be the Coxeter diagram  of $W$. 
We take $Y$ to be an affine subdiagram of $X$, $p$ to be a node in
$X\backslash Y$, and repeat the argument that  proves proposition \ref{sbetas} above. \qed
\begin{remark}
1. In view of remark \ref{light}, corollary \ref{expo} 
also holds for non simply laced, irreducible Coxeter groups $W$ whose
Coxeter diagrams contain a proper affine subdiagram. 

\noindent 2. It is easy to see that ({\em cf \S \ref{w3prop-sec} below}) that $\wth$
has exponential growth. Thus this proposition gives another proof (in
the simply laced case) of the
result mentioned in the introduction: {\em an irreducible Coxeter
  group $W$ which is infinite and non-affine has exponential growth}.

\end{remark}

\section{Properties of $\wth$}\label{w3prop-sec}
To complete the proof of proposition \ref{mainprop}, we must study the
reflection subgroup constructed in proposition \ref{sbetas} more closely.
We collect together some useful properties of the Coxeter group
$\wth$. Note that $\wth$ is just the free product of three groups of
order 2. The following facts are all fairly standard, and we omit proofs:
\begin{proposition}\label{conjaction}
\be
\item The Poincar\'{e} series $\wth(q) = \dfrac{1+q}{1-2q} $.
\item \label{one} Each $w \in \wth$ has a unique reduced expression as
  a product of Coxeter generators.
\item If $\sth:=\{s_1,s_2,s_3\}$ are the Coxeter generators of $\wth$, the
  conjugacy classes of the $s_i$  are pairwise disjoint. Further $w
  \in \wth$, $w s_i w^{-1} = s_i \Leftrightarrow w \in \{1,s_i\}$.
\ee
\end{proposition}

We let $\tth:=\bigcup_{w \in \wth} w \sth w^{-1}$ be the set of
reflections in $\wth$.
Proposition \ref{conjaction} implies that $\tth$ is a disjoint union of
the orbits of the $s_i,\, (i=1,2,3)$ under the conjugation action of
$\wth$; further, the stabilizer of $s_1$ is $\{1, s_1\}$.
 Let $\oh \subset \tth$ be the orbit of
$s_1$; if we let $K:=\{s_1\} \subset \sth$, then $\oh = \{\sigma s_1
 \sigma^{-1}: \sigma \in \wth^K\}$, with $ \sigma s_1 \sigma^{-1} \neq
 \tau s_1 \tau^{-1}$ for $\sigma \neq \tau \in \wth^K$; here $\wth^K$
 is the set of minimal left coset representatives of the parabolic subgroup
 $(\wth)_K =\{1,s_1\}$.

The Poincar\'{e} series of $\wth^K$ is $\wth(q)/(1+q) =
1/(1-2q)$; so for each $k \geq 0$, there are $2^k$ elements $\sigma \in \wth^K$ such that
$\ell(\sigma)=k$. For these $\sigma$, $\ell(\sigma s_1 \sigma^{-1})
\leq 2\ell(\sigma) + 1 = 2k+1$. So, 
\bprop \label{twok}
For each $k \geq 0$, there are $\geq 2^k$ elements $t \in \oh$ such that $\ell(t) \leq 2k+1$.
\eprop
Let $\vth$ be the geometric representation of $\wth$ with basis $\{\alpha_1, \alpha_2, \alpha_3\}$ and 
invariant bilinear form $(,)$. We remark that there are many choices
for the $\wth$ invariant form $(,)$ on $\vth$. It only needs to
satisfy $(\alpha_i, \alpha_i) = 1 \,\forall i$ and $(\alpha_i, \alpha_j) \in \integers^{\leq -1} \, i \neq j$.
Let $\rs(\wth) \subset \vth$ be the root system of $\wth$.  We then have:
\bprop \label{nonzerocoeff}
Let $\alpha \in \rs[+](\wth)$ such that $s_{\alpha} \in \oh$. If
$\alpha = \sum_{i=1}^3 c_i \alpha_i$ ($c_i \in \nat$), then $c_1 >0$.
\eprop
\noindent
\pfbegin Given $\gamma_1, \gamma_2 \in \rs(\wth)$ write $\gamma_1 >
\gamma_2$ if $\gamma_1 - \gamma_2$ is a nonnegative integer linear
combination of the $\alpha_i$. It is a well known fact
(see for e.g the argument used in \cite[proposition 5.1(e)]{kac})
that given a positive root $\alpha$, there exists a sequence $\gamma_0 >
\gamma_1 > \cdots > \gamma_r$ such that (i) $\gamma_0 = \alpha$,
$\gamma_r \in \{\alpha_i:i=1,2,3\}, \,\gamma_j \in \rs[+](\wth) \forall
j$ (ii) For each $p$, $\gamma_{p+1} = s_{i_p} (\gamma_p) \text{ for
  some } i_p\in\{1,2,3\}$. 

Thus each $\gamma_p \in \wth.\alpha$ or equivalently $s_{\gamma_p}$ is
$\wth$ conjugate to $s_\alpha$. The disjointness of the orbits of the
$s_i$ mentioned before and the hypothesis that $s_\alpha  \in \oh$ imply that $\gamma_r = \alpha_1$.
So $\alpha = \alpha_1 + \sum_{p=0}^{r-1} (\gamma_p -
\gamma_{p+1})$. \qed

\section{Proof of proposition \ref{mainprop}}
We now put together the results of the previous two sections. Let $W,
S, J, \beta_i$ be as in \S \ref{thebetas}. Let $W\p
= \langle  s_{\beta_i}: i=1,2,3 \rangle$ be the reflection subgroup isomorphic to
$\wth$ constructed in Proposition \ref{sbetas}. Let $S\p
:=\{s_{\beta_i}: i=1,2,3\}$. Define $\rs(W\p) := W\p.\{\beta_i\}_{i=1}^3
\subset \rs(W)$. We identify $\rs(\wth)$ with $\rs(W\p)$ by sending
$\alpha_i \mapsto \beta_i$ and requiring that this map commute with
the $\wth$ action (for this identification to be a linear map of the
underlying vector spaces, we will need to use the form on $\vth$
that satisfies $(\alpha_i, \alpha_j) = (\beta_i, \beta_j) \; \forall
i,j$). Now, applying propositions  \ref{twok} and
 \ref{nonzerocoeff} to $W\p \cong \wth$, we deduce that for each $k \geq 0$, there are $\geq 2^k$
elements $\beta \in \rs[+](W\p) \subset \rs[+](W)$ such that
$\ell_{S\p}(s_{\beta}) \leq 2k+1$ and $\beta = \sum_{i=1}^3 c_i
\beta_i$ with $c_1 > 0$. Here $\ell_{S\p}(\cdot)$ denotes the length
function of $W\p$ w.r.t $S\p$. Now let $M:=\max\{\ell(s_{\beta_i}):
i=1,2,3\}$, where $\ell(\cdot)$ is the usual length function on $W$
w.r.t $S$. We clearly have $\ell(w) \leq M \ell_{S\p}(w) \; \forall w
\in W\p$. Thus for the $\beta$'s above, $\ell(s_{\beta}) \leq M(2k+1)$.
We now {\em Claim:} for each of the above $\beta$'s, $s_{\beta} \in T
\backslash W_J$.
Referring back to the statement of Proposition~ \ref{mainprop}, we see
that this claim together with  what we have shown thus far would complete
the proof of that proposition.

{\em Proof of Claim:} Recall from \S \ref{thebetas} that $p$ was
chosen to be an element of $S\backslash J$, and that $\beta_1 = \sum_{q \in S}
d_q \alpha_q$ with $d_p = 1$. Now each of the $\beta$'s of the above
paragraph can be written as $\beta=\sum_{i=1}^3 c_i \beta_i$ with
$c_1 >0$. It follows then that we can write $\beta = \sum_{q \in S}
e_q \alpha_q$ with $e_p > 0$. Since $p \not\in J$, this means that
$\beta$ is not a linear combination of the simple roots $\alpha_q, q
\in J$. It is an easy fact that this implies $s_\beta \not\in W_J$
(sketch of proof: If $ w \in W_J, w(\beta) = \beta - (\text{a linear
combination of } \alpha_q, q \in J) = \sum_{u \in S} k_u \alpha_u$
  with $k_p = e_p >0$. Thus $w(\beta) \in \rs[+](W),  \forall w \in
  W_J$. But $s_\beta(\beta) = -\beta \in \rs[-](W)$; this gives the
  desired contradiction).

Thus, putting everything back together,  our main theorem \ref{mainthm} is proved. \qed

\section{Quotients by reflection subgroups}
We assume $(W,S)$ to be a simply laced, irreducible Coxeter system
which is neither finite nor affine. As usual, we let $V$ be its
geometrical realization, $(\cdot,\cdot)$ the invariant bilinear form etc. Let
$W\p$ be a finitely generated reflection subgroup of $W$ and $S\p =
\{s_{\beta_i}\}_{i=1}^k$ be its Coxeter generators as in theorem
\ref{deo-dyer}, with $\beta_i \in \rs[+](W)$. Let
$\rs(W\p):=W\p.\{\beta_i\}_{i=1}^k$ be its root system. 

It was shown by Dyer \cite[(3.4)]{dyer1} that the left cosets of $W\p$
in $W$ have {\em unique} elements of minimal length; these elements
$w$ are determined by the condition that $\ell(ws_{\beta_i}) >
\ell(w)\; \forall i$ or equivalently by the condition $w(\beta_i) \in
\rs[+](W) \;\forall i$. We let $\lc$ denote the set of these minimal
coset representatives. We remark that while each $w \in W$ can be
uniquely written as $w=\sigma \tau$ with $\sigma \in \lc, \tau \in
W\p$, it may no longer be true that $\ell(w) = \ell(\sigma) +
\ell(\tau)$. The natural question now is to study the growth of
$W/W\p$ or more precisely, the growth of $\lc$.  

To make our arguments simpler, we assume further that $W\p$ is
irreducible as a Coxeter group, leaving the details of the reducible
case to the reader. We first consider an additional hypothesis on the
$W\p$:
\begin{lemma}\label{class}
Let notation be as above. TFAE:
\be
\item \label{clone} $\exists \,\alpha \in \rs[+](W)\backslash
  \rs[+](W\p)$ such that its $W\p$ orbit satisfies $\{\alpha\}
  \subsetneq W\p \cdot \alpha \subset \rs[+](W)$.
\item \label{cltwo} $\exists \,\gamma \in \rs[+](W)\backslash
  \rs[+](W\p)$ such that $(\gamma, \beta_i) \leq 0 \; \forall
  i=1\cdots k$, with  not all $ (\gamma, \beta_i)$ equal to zero.
\ee
\end{lemma}

\noindent
{\em Proof:} $\eqref{clone} \Rightarrow \eqref{cltwo}$: For a positive
root $\beta = \sum_{s \in S} c_s \alpha_s \in \rs[+](W)$, define its
{\em height} to be  $\sum c_s (> 0)$. Given $\alpha$ as in \eqref{clone}, let $\gamma$ be an
element in the $W\p$ orbit $W\p \cdot \alpha$ of minimal height. It is
clear that $\gamma$ satisfies $\eqref{cltwo}$. For \eqref{cltwo}
$\Rightarrow$ \eqref{clone}: given such $\gamma$, it is a standard
argument that $\forall \sigma \in W\p$, $\gamma - \sigma \gamma =
\sum_{i=1}^k d_i \beta_i$ with $d_i \leq 0 \; \forall i$. Thus
$\sigma \gamma \in \rs[+](W) \; \forall \sigma \in W\p$. So $\alpha :=
\gamma$ satisfies $\eqref{clone}$. \qed

\vspace{0.2cm}
The next proposition explains the usefulness of the above conditions:
\bprop
If $W\p$ satisfies the equivalent conditions of lemma \ref{class},
then $\lc$ has exponential growth.
\eprop

\noindent
{\em Proof:}
First observe that if $W\p$ is finite or affine, the truth of the proposition
follows from the exponential growth of $W$ coupled with the polynomial
growth of $W\p$. So, assume $W\p$ is neither finite nor affine.

Let $\{s_{\beta_i}\}_{i=1}^k$ be the Coxeter generators of $W\p$. 
Pick $\gamma$ as in condition \eqref{cltwo} of the lemma and define
$\ovw$ to be the reflection subgroup generated by $s_{\gamma}$ and
$W\p$. By Dyer's criterion (theorem \ref{dyer}), 
$\{s_{\gamma}, \,s_{\beta_i} (i=1\cdots k)\}$ are
the Coxeter generators of $\ovw$. Now, $\ovw$ is an irreducible
Coxeter group which is neither finite nor affine and this contains
 $W\p$ as a parabolic subgroup. By our main theorem \ref{mainthm}, the
set of minimal coset representatives of $W\p$ in $\ovw$ has
exponential growth (wrt the length function on $\ovw$). Thus, if
$$a_m:=\# \{\sigma \in \lc \cap \ovw: \ell_\ovw(\sigma) \leq m\}$$
then $\limsup_{m \to \infty} a_m^{1/m} > 1$. If
$K:=\max\{\ell(s_\gamma), \,\ell(s_{\beta_i}) (i=1\cdots k)\}$, we have
$\ell(\sigma) \leq K \ell_{\ovw} (\sigma) \; \forall \sigma \in
\ovw$. So, if $b_m:=\# \{\sigma \in \lc: \ell(\sigma) \leq m\}$, then
$b_{Km} \geq a_m$. Thus
$$\limsup_{m \to \infty} b_m^{\frac{1}{m}}
 \geq \limsup_{m \to \infty} \left( b_{Km} \right) ^{\,1/Km}
 \geq  \left( \limsup_{m \to \infty} a_m^{\frac{1}{m}}\right)^{\frac{1}{K}} > 1$$
Thus $b_m$, and hence $\lc$, has exponential growth. \qed

\vspace{0.2cm}
\noindent
{\bf Final Remarks:}
It is easily seen if  $W\p = W_J \subsetneq W$ or $W\p =\sigma W_J
 \sigma^{-1}$ (parabolic subgroups and their conjugates), then the
 equivalent conditions of lemma \ref{class} are satisfied. In
 practice, when $W$ has small rank and given the Coxeter generators
 $s_{\beta_i}$ of $W\p$ explicitly, it is often easy to show that this lemma is
 satisfied by explicitly producing a positive root $\gamma$ such that $(\gamma,
 \beta_i) \leq 0 \, \forall i$. From such examples worked out
 by hand, it appears that this lemma is satisfied for a large number of
 reflection subgroups of $W$ (when $W$ is non-finite, non-affine).
It would useful to be able to completely characterize such reflection
subgroups.

\bibliographystyle{amsplain}
\bibliography{biblio}
\end{document}